# Topology optimization of periodic lattice structures for specified mechanical properties using machine learning considering member connectivity

Tomoya Matsuoka[1], Makoto Ohsaki[1] and Kazuki Hayashi[1]

[1] *Department of Architecture and Architectural Engineering, Kyoto University, Kyoto 615-8540, Japan*

**Abstract**
This study proposes a methodology to utilize machine learning (ML) for topology optimization of periodic lattice structures. In particular, we investigate data representation of lattice structures used as input data for ML models to improve the performance of the models, focusing on the filtering process and feature selection. We use the filtering technique to explicitly consider the connectivity of lattice members and perform feature selection to reduce the input data size. In addition, we propose a convolution approach to apply pre-trained models for small structures to structures of larger sizes. The computational cost for obtaining optimal topologies by a heuristic method is reduced by incorporating the prediction of the trained ML model into the optimization process. In the numerical examples, a response prediction model is constructed for a lattice structure of 4 × 4 units, and topology optimization of 4 × 4-unit and 8 × 8-unit structures is performed by simulated annealing assisted by the trained ML model. The example demonstrates that ML models perform higher accuracy by using the filtered data as input than by solely using the data representing the existence of each member. It is also demonstrated that a small-scale prediction model can be constructed with sufficient accuracy by feature selection. Additionally, the proposed method can find the optimal structure in less computation time than the pure simulated annealing.

**Keywords:** *Lattice structure, topology optimization, machine learning, filtering process, feature selection*

## 1. Introduction

Lattice structures considered in this study are planar periodic frame structures composed of squares with diagonals as a basic pattern. Various mechanical properties, such as elastic modulus, can be achieved by optimizing the topology of lattice structures so as to have the target deformation against specific loads or the target reaction forces for specified displacements [1, 2]. Optimization of lattice structures can be considered either on micro scales, as in materials design, or on macro scales, as in the design of building structures. Lattice structures considered in this study are intended for the macro scales and can serve as components for structural walls or reinforcement. The optimized structures consist of lightweight structural units with a given performance.

In this study, among several structural optimization types, we focus on topology optimization, in which the arrangement of structural elements is to be designed. Since topology optimization has discrete variables representing the presence or absence of each structural element, it is difficult to use an optimization method based on gradients, which necessitates a large number of functional evaluations. For this reason, research on replacing a part of the optimization process with machine learning (ML) models have become popular in recent years.

The mechanical properties of structures are determined not only by the presence or absence of each structural element but also by the relationship with adjacent structural elements, and many methods have been proposed to consider such relationships in the construction of response prediction models. One such method is convolutional neural networks (CNNs). Yang et al. [3] predicted stress-strain curves of binary composites, and Takagi et al. [4] predicted mechanical properties of lattice structures whose structural patterns can be represented by a two-dimensional array. However, in general, lattice structures cannot represent their structural patterns in a two-dimensional array, making it difficult to use CNNs such as those used in image recognition.

The filtering method processes the arrangement of multiple members, which allows for considering the connectivity of general braced frame structures [5]. In this method, convolution operations transform data into a combinatorial representation of the existence of multiple members. Since the filtered data has larger dimensions than the data in which the absence or presence of each member is represented by 0-1, it increases the computational cost for training the ML model and may even worsen the prediction accuracy. Therefore, feature selection plays a significant role in determining particularly important local patterns of member arrangement for the ML models.

In this study, we propose a method for constructing a response prediction neural network (NN) model for periodic lattice structures, using data generated by the filtering process and feature selection. We also present a convolution approach to apply pre-trained models for small structures to structures of larger sizes by converting the data for a large structure into the data representing an equivalent small structure. Furthermore, we demonstrate topology optimization of periodic lattice structures using a heuristic method that incorporates NN-based identification of good solutions.

This paper is organized as follows. In Sec. 2 we present the periodic lattice structures and the topology optimization problem with the compliance to horizontal forces as the objective function. In Sec. 3 a method is proposed for constructing NNs to predict the responses of the lattice structures, focusing on the filtering process and feature selection. In Sec. 4 we describe a heuristic method used for solving the topology optimization problem, and the procedure to incorporate the prediction results of the trained NN models. In Sec. 5 response prediction NN models are constructed for $4 \times 4$-grid units and $8 \times 8$-grid units, examining the effects of the filtering process and feature selection, and their frame arrangements are then optimized as numerical examples. Sec. 6 concludes the paper.

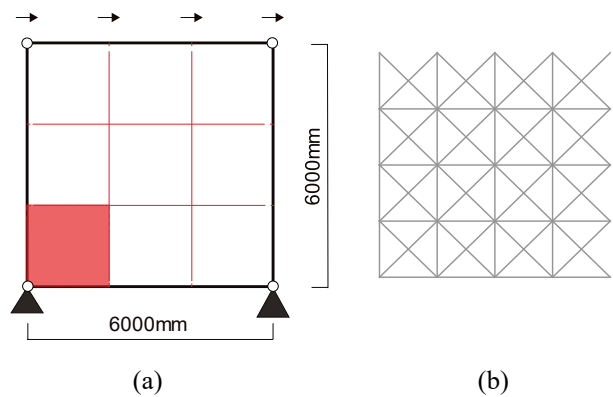

(a)          (b)

Figure 1. A lattice structure model; (a) whole structure, (b) ground structure of $4 \times 4$-grid unit

## 2. Topology optimization of periodic lattice structures

Consider a structure with 9 identical units subjected to horizontally distributed loads as shown in Fig. 1(a). The structure consists of the unit of $m \times m$ grids as the ground structure. The ground structure of the $4 \times 4$-grid unit is shown in Fig. 1(b).

The conditions for the structural analysis are as follows:

- Frame members are modeled by Euler-Bernoulli beam elements.
- Rigid members are placed on the outer frame of the structure.
- Only the nodes at the four corners of the model are pin-connected, and all other nodes are rigidly connected.
- Diagonal members are not connected at the center of the grid.
- A horizontal load of 12 kN in total is distributed to the nodes on the upper edge.
- The cross-section and material properties of bar members are determined as Table 1.

Table 1. Specification of beam elements

|  | Width $B$ (mm) | Depth $D$ (mm) | Cross-sectional area $A$ (mm²) | Second moment of area $I$ (mm⁴) | Young's modulus $E$ (N/mm²) |
|---|---|---|---|---|---|
| Outer frame | 300 | 400 | $1.2 \times 10^5$ | $1.6 \times 10^9$ | 20000 |
| Lattice frame ($m \times m$ grid) | 150 | $400/m$ | $6.0 \times 10^4/m$ | $8.0 \times 10^8/m^3$ | 20000 |

We consider the following topology optimization problem, where the variable $x_i$ represents the presence or absence of member $i$, the upper bound of the total volume of the unit $v(x)$ is $v_c$, and the compliance of the whole structure to horizontal forces is the objective function to be minimized.

$$\text{Minimize } \boldsymbol{u}(\boldsymbol{x})^\top \boldsymbol{F} \tag{1a}$$
$$\text{subject to } x_i \in \{0,1\} \tag{1b}$$
$$v(\boldsymbol{x}) \leq v_c \tag{1c}$$

## 3. Neural networks for response prediction of lattice structures

### 3.1. Generation of training data

We randomly generate 150000 combinations of 0-1 variables representing the absence or presence of members in a $m \times m$-grid unit. The corresponding structures are then analyzed to determine their compliance to horizontal forces. Structures with a compliance greater than a threshold value of 5 kN are excluded from the data as unstable.

### 3.2. Filtering process

Although existence of bar members can be described by 0-1 variables, this binary representation cannot consider the connectivity of multiple members. Therefore, ML models may not accurately capture the mechanical properties of lattice structures. To explicitly consider the combinations of adjacent members, we use the data processed by filtering as an input for the ML models.

The combination of multiple members is represented regarding a subregion consisting of 8 members connected to each node of the structure as shown in Fig. 2. The number of subregions is denoted as $n_\mathrm{n}$, and the combination of $n_\mathrm{m}$ ($< 8$) members are considered. A filter matrix $\mathbf{C} \in \mathbb{R}^{{}_8C_{n_\mathrm{m}} \times 8}$ is generated to store all combinations with $n_\mathrm{m}$ elements valued 1 and $8 - n_\mathrm{m}$ elements valued 0, along with a matrix $\mathbf{X}_0 \in \mathbb{R}^{n_\mathrm{n} \times 8}$ that stores the presence or absence of members for each subregion. Then, filtering is performed to convert $\mathbf{X}_0$ to $\mathbf{X}$ as follows:

$$\mathbf{X} = \lfloor \mathbf{X}_0 \mathbf{C}^\top / n_\mathrm{m} \rfloor \tag{2}$$

where $\lfloor \cdot \rfloor$ is the floor operator to have the largest integer that does not exceed the real value.

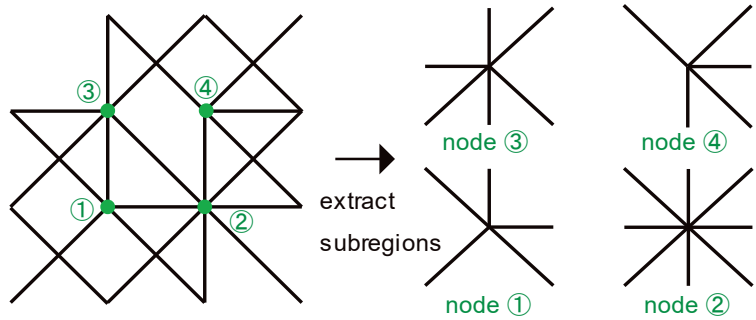

Figure 2. Extracting subregions from lattice structures for filtering

### 3.2. Feature selection

The filtering process described in Sec. 3.1 explicitly represents all combinations for a specified number of members, but not all information is necessarily required. Building ML models using data with such large dimension leads to higher computational cost due to a large number of parameters to be updated and may even worsen the prediction accuracy. Therefore, feature selection is performed to select important features based on statistical evaluation.

In this study, feature selection is performed using the f_regression function of scikit-learn [6], a function that statistically evaluates the linear dependence between each feature and the objective function value. Suppose the number of data is $n$, the total number of features is $k$, and the correlation coefficient between each feature $X_i$ and the objective function $y$ is $r_{X_i y}$. Then, the F-statistic is calculated as follows:

$$F_{X_i} = \frac{r_{X_i Y}^2}{1 - r_{X_i Y}^2} \cdot \frac{n - k - 1}{k} \tag{3}$$

The larger the value of the F-statistic, the greater the correlation between each feature and the objective function value. It is important to identify the minimum number of features in the descending order of their F-statistic value so that the prediction error of the ML model is not significantly impaired.

### 3.3. Neural networks

In this study, an NN model is used as the regression model and implemented in PyTorch [7].

The NN consists of 3 intermediate fully-connected layers with units 900-600-300, and the activation function is the identity function in the output layer and the Rectifier Linear Units (ReLU) in the other layers. For the training condition, the minibatch size is 200, the number of epochs is 100, the error function is the mean squared error (MSE), and the parameter optimization algorithm is Adam [8].

In the NN construction method described above, the trained NN model corresponds only to structures of a specific size. Therefore, we propose the method to convert the data representing a finer lattice pattern into the data representing an equivalent coarse lattice pattern by convolution operation. Then, the pre-trained model for small structures can be applied to structures of larger sizes.

The members of the lattice structures are classified into 4 types according to their orientation: horizontal, vertical, right diagonal and left diagonal. By independently expressing the existence of the members across 4 channels, we can perform a two-dimensional convolution operation, since the members of the same channel are regularly arranged in a two-dimensional array, as illustrated in Fig. 3. A $2m \times 2m$-grid unit represented by a three-dimensional array of sizes $(4, 2m, 2m)$ is artificially transformed into an $m \times m$-grid unit represented by an array of sizes $(4, m, m)$ through convolution. The convolution is processed from 4-channel inputs to 4-channel outputs with a kernel of size $2 \times 2$. This output array is vectorized and can be used as input data for the NN model pre-trained for lattice structures with $m \times m$ grids in the unit.

The artificially converted $m \times m$-grid unit has similar structural characteristics to the actual $m \times m$-grid unit, allowing for the same filtering and feature selection processes described above. However, the filtering process in Sec. 3.2 assumes binary variables, whereas the transformed data after convolution contains real values. As a result, the filter processing expressed by Eq. (2) might miss information with regard to the arrangement of members unless the two-dimensional array is filled with 1. Therefore, we use a Sigmoid-like function (Eq. (4)) to treat data continuously, and we perform the filtering process as Eq. (5). Figure 4 shows the whole network structure.

$$f(x; n_\mathrm{m}) = \frac{1}{1 + \exp(-10(x - 0.75 n_\mathrm{m}))} \tag{4}$$

$$\mathbf{X} = f(\mathbf{X}_0 \mathbf{C}^\mathsf{T}; n_\mathrm{m}) \tag{5}$$

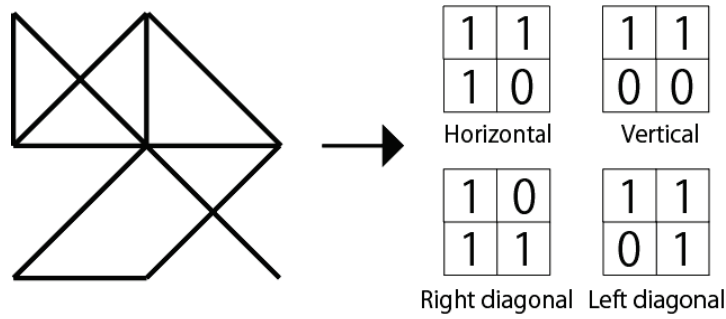

Figure 3. Representation of member arrangement using 4 channels

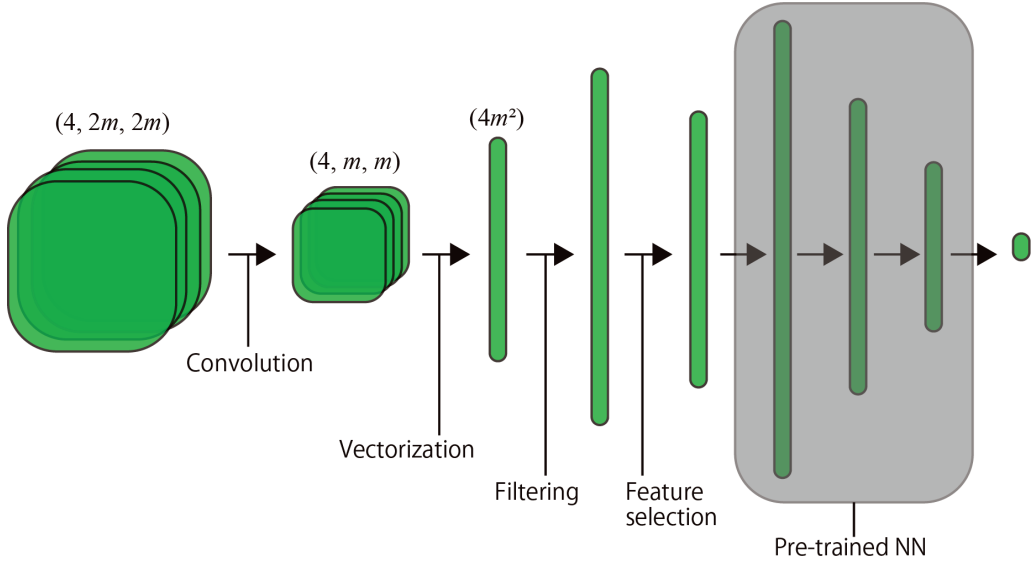

Figure 4. Whole network structure applied to the pre-trained NN

## 4. Heuristic optimization method

The optimal solutions are searched using simulated annealing (SA), which is a heuristic method that stochastically updates a solution. While SA is a highly versatile method that can be applied to various optimization problems, it requires an enormous number of function evaluations. Therefore, we propose an algorithm that omits structural analysis for inferior solutions in the iterations based on the prediction by the NN model.

In addition, in order to solve the optimization problem with constraints in Eq. (1) with SA, we consider the following objective function that has the penalty term $\lambda(v(\pmb{x}) - v_c)$ with a positive penalty coefficient $\lambda$ for the solutions that violate the constraints.

$$g(\pmb{x}; \lambda) = \pmb{u}(\pmb{x})^\top \pmb{F} + \lambda \max(v(\pmb{x}) - v_c, 0) \tag{6}$$

The SA algorithm assisted by the trained NN model for the minimization of penalized objective function $g(\pmb{x}; \lambda)$ is described as follows:

| SA algorithm assisted by the trained NN model |
|---|
| Input: $x_0$ (initial solution), $a$ (cooling coefficient), $n_s$ (number of cooling steps), $p$ (cooling period), $n_v$ (maximum number of variables to change), $w_d$ (decrement of threshold value), $w_i$ (increment of threshold value) |
| Output: optimal solution |

| | |
|---|---|
| Step1. | Calculate the initial objective function value as $y_0 = g(x_0; \lambda)$, and set an initial threshold value for predictions that determine necessity of structural analysis as $c = y_0$. Set a scaling parameter as $s = -0.1y/(\log 0.5)$, and the initial temperature as $T_0 = 1.0$. Initialize the iteration counter as $k = 0$. |
| Step2. | Increase the counter as $k \leftarrow k + 1$, and set $x_k = x_{k-1}, y_k = y_{k-1}$. Randomly choose a subregion, and switch between 0 and 1 for at most $n_v$ variables in the subregion to generate the neighborhood solution $\tilde{x}$. |
| Step3. | Predict the value of objective function $\tilde{y}_p$ by using the NN model for $\tilde{x}$. If $\tilde{y}_p < c$, decrease $c$ to $c - w_d$ and go to Step 4; otherwise, increase $c$ to $c + w_i$ and go to Step 7. |
| Step4. | Calculate the objective function value $\tilde{y} = g(\tilde{x}; \lambda)$ by carrying out structural analysis. |
| Step5. | Calculate the transition probability to the neighborhood solution $\tilde{x}$ by following equation. $$p = \begin{cases} 1 & (\tilde{y} < y_{k-1}) \\ \exp\left(\frac{y_{k-1} - \tilde{y}}{Ts}\right) & (\tilde{y} \geq y_{k-1}) \end{cases} \quad (7)$$ |
| Step6. | Generate a random number $r \in [0,1)$ and set $x_k = \tilde{x}, y_k = \tilde{y}$ if $r < p$. |
| Step7. | Update the temperature $T$ to $aT$ after $p$ iterations at the same $T$. Return to Step 2 when $k < pn_s$; otherwise, terminate the algorithm and output the best solution as the optimal solution. |

## 5. Numerical examples

### 5.1. Construction of response prediction neural networks

#### 5.1.1. Generation of training data

The datasets for $4 \times 4$-grid units and $8 \times 8$-grid units are generated by the method described in Sec. 3.1. Figure 5 shows the relationship between the compliance and structural volume of the unit as scatter plots, and Table 2 shows the statistics for compliance and volume in the generated datasets in comparison with the ground structure.

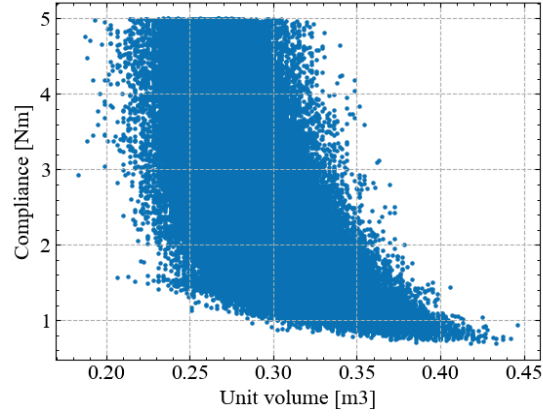

(a)

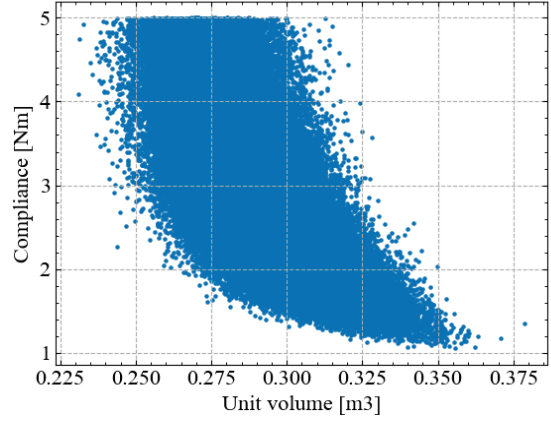

(b)

Figure 5. Relationship between volume and compliance; (a) 4 × 4-grid unit, (b) 8 × 8-grid unit

Table 2. Statistics for volume and compliance in the datasets along with the values of the ground structure

|  |  | Average of the data | Standard deviation of the data | Ground structure |
|---|---|---|---|---|
| 4 × 4-grid unit | Volume (m³) | 0.29931 | 0.03130 | 0.57941 |
| (123168 samples) | Compliance (Nm) | 2.29259 | 0.99183 | 0.53748 |
| 8 × 8-grid unit | Volume (m³) | 0.29235 | 0.01671 | 0.57941 |
| (135364 samples) | Compliance (Nm) | 2.77769 | 0.84768 | 0.54031 |

5.1.2. Neural networks for 4 × 4-grid unit

We examine the effect of the filtering process described in Sec. 3.2 on the construction of the NN model. 75% of the data in Sec. 5.1.1 is used as the training data and 25% as the test data. Figure 6 shows the history of MSE loss in the learning process for each type of filtering. In Fig. 6, "No filtering" refers to the case where the data simply indicates the presence or absence of each member, and "Filtering 2" and "Filtering 3" refer to the case where the data is filtered as $n_\mathrm{m} = 2$ and $n_\mathrm{m} = 3$, respectively. From the figure, the case $n_\mathrm{m} = 2$, which considers the combination of two members, leads to the best prediction accuracy. Therefore, in the examples below, filtering is processed with $n_\mathrm{m} = 2$.

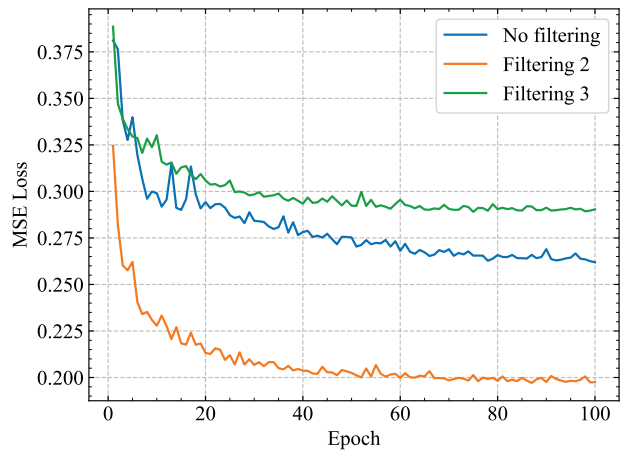

Figure 6. Comparison of accuracy of NN model with each filtering type

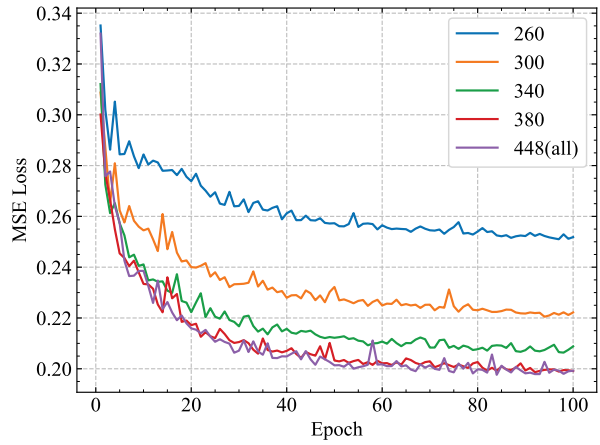

Figure 7. Comparison of accuracy of NN model with five different numbers of input features

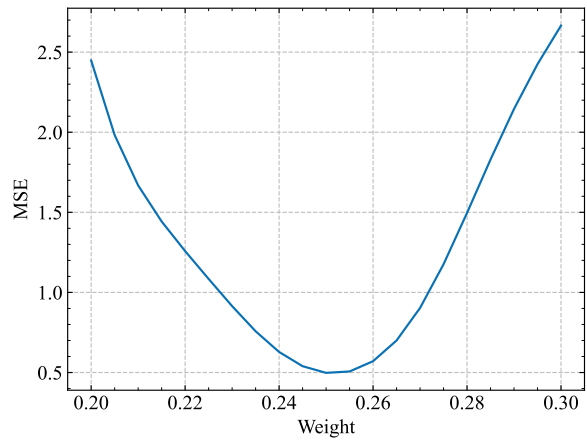

Figure 8. Prediction error for 8 × 8-grid units with respect to each weight parameter

The effect of feature selection described in Sec. 3.3 on the filtered data is then examined. The size of filtered input data is 448 because 16 subregions and $_8C_2$ combinations of bar members are considered. Figure 7 shows the training history of the NN model using 5 types of input data with different numbers of selected features. It is observed from the figure that the accuracy declines as the number of features decreases; however, the decline in accuracy is relatively small if the 340 or 380 features are used. Therefore, the data with 340 features is used hereafter.

5.1.3. Neural networks for 8 × 8-grid unit

The NN model for 4 × 4-grid units, built using the filtering process and feature selection in Sec. 5.1.2, is used to predict the response of 8 × 8-grid units as described in Sec. 3.3. This method does not require training data because it is not necessary to update the parameters of NN, but the MSE is measured with the data generated in Sec. 5.1.1 for verification purpose. Figure 8 shows the prediction error for each weight parameter in the convolution. The figure shows that the best accuracy is achieved at a weight of 0.25, which is reasonable because the convolution is carried out on 4 elements and specification of beam elements are determined in Table 1 so that 4 × 4-grid and 8 × 8-grid units are mechanically equivalent.

Table 3. Results of topology optimization of 4 × 4-grid unit

| Initial solution | SA assisted by the trained NN model | | Pure SA | |
|---|---|---|---|---|
| | Compliance (Nm) | Volume (m³) | Compliance (Nm) | Volume (m³) |
| 1. | 0.79654 | 0.26024 | 0.80036 | 0.26024 |
| | 0.72302 | 0.25456 | 0.72246 | 0.25456 |
| | 0.79317 | 0.26024 | 0.78388 | 0.25585 |
| | <span style="color:red">0.71656</span> | <span style="color:red">0.25456</span> | 0.72261 | 0.26024 |
| | 0.79411 | 0.26023 | 0.78531 | 0.26024 |
| 2. | 0.79573 | 0.26024 | 0.71630 | 0.25456 |
| | 0.78665 | 0.25585 | 0.78867 | 0.25585 |
| | 0.80094 | 0.26024 | 0.80878 | 0.25585 |
| | 0.78594 | 0.26024 | 0.79992 | 0.26024 |
| | 0.79676 | 0.25713 | 0.78604 | 0.25585 |
| 3. | 0.80535 | 0.26024 | 0.72302 | 0.25456 |
| | 0.80422 | 0.26024 | 0.79777 | 0.26024 |
| | 0.80212 | 0.26024 | <span style="color:red">0.71157</span> | <span style="color:red">0.25456</span> |
| | 0.80576 | 0.25713 | 0.80495 | 0.26024 |
| | 0.80344 | 0.25895 | 0.80530 | 0.26024 |
| 4. | 0.80298 | 0.25585 | 0.78830 | 0.26024 |
| | 0.79083 | 0.25585 | 0.80769 | 0.25713 |
| | 0.80836 | 0.26024 | 0.78876 | 0.26024 |
| | 0.80576 | 0.26024 | 0.78594 | 0.26024 |
| | 0.79107 | 0.25584 | 0.79004 | 0.25585 |
| Average | 0.79047 | 0.25842 | 0.77588 | 0.25785 |
| Average number of structural analyses | 8853 | | 64000 | |

## 5.2. Topology optimization

### 5.2.1. 4 × 4-grid unit

Topology optimization of a 4 × 4-grid unit with $v_c = 0.26$ is performed by the SA algorithm with and without the trained NN model. The parameters are set to $\lambda = 10, a = 0.88, n_s = 50, p = 640, n_v = 3, w_d = 0.001$, and $w_i = 0.0003$. Table 3 shows the optimization results for 4 different initial solutions, with each undergoing optimization 5 times. Among the obtained solutions, the solution that strictly meets the volume constraint and achieves the lowest compliance is highlighted in red. The best structure obtained by the SA with NN model is shown in Fig. 9.

According to the Table, the SA algorithm assisted by the trained NN model achieves the optimization results comparable to the pure SA, while omitting a great number of structural analyses.

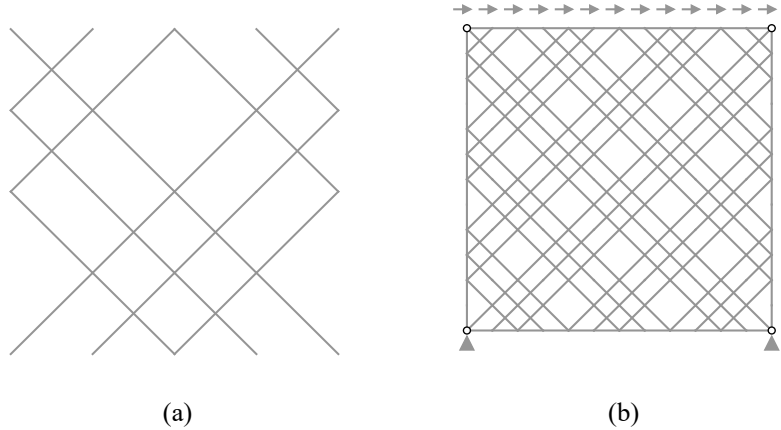

(a)                         (b)

Figure 9. Optimal lattice structure of 4 × 4-grid unit; (a) unit (b) whole structure

### 5.2.2. 8 × 8-grid unit

It is more challenging to optimize the entire structure with the 8 × 8-grid units, as it contains 256 members in each unit. Therefore, we utilize the optimal structure of the 4 × 4-grid unit and search for the corresponding near-optimal structure. First, an 8 × 8-grid initial structure is generated by locally duplicating each member of the 4 × 4-grid structure into 4 members. A local search is then performed. The local search does not allow transitions to a solution with a worse objective function value, which is equivalent to running the SA algorithm with $T = 0$. The parameters are set to $\lambda = 10, p = 12800, n_v = 3, w_d = 0.001$, and $w_i = 0.0003$. Figure 10 shows the structure of initial solution with the volume of 0.25456m³ and the compliance of 0.81299Nm corresponding to the optimal structure in Fig. 9. The results of 5 optimization trials are shown in Table 4, where the feasible solution with the lowest compliance is highlighted in red. The best structure obtained by the local search with NN model is shown in Fig. 11.

Optimization results are not much different with and without the NN model, while omitting a great number of structural analyses. Comparing Figs. 9 and 11, we can see that they have similar lattice patterns and that the local search described above is effective for finding a near optimal solution.

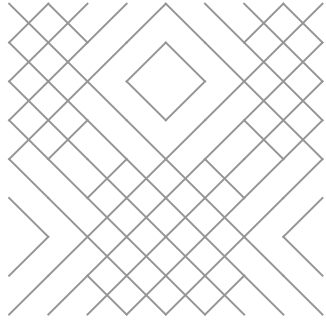

Figure 10. Initial structure for local search of 8 × 8-grid unit corresponding to the best structure of 4 × 4-grid unit

Table 4. Results of local search for near-optimal structure of 8 × 8-grid unit

|  | SA assisted by the trained NN model | | Pure SA | |
| --- | --- | --- | --- | --- |
|  | Compliance (Nm) | Volume(m³) | Compliance (Nm) | Volume(m³) |
|  | 0.79564 | 0.25992 | 0.78829 | 0.25960 |
|  | 0.80116 | 0.25960 | 0.79275 | 0.25992 |
|  | 0.78582 | 0.25960 | 0.78653 | 0.25960 |
|  | 0.78918 | 0.25927 | 0.78126 | 0.25960 |
|  | 0.79409 | 0.25927 | 0.78611 | 0.25992 |
| Average | 0.79318 | 0.25953 | 0.78699 | 0.25973 |
| Average number of structural analyses | 2835 | | 12800 | |

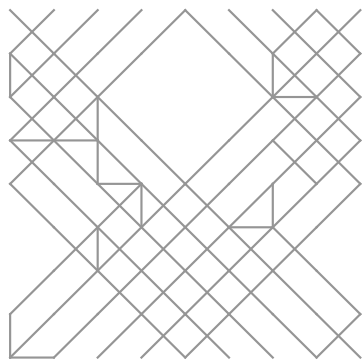 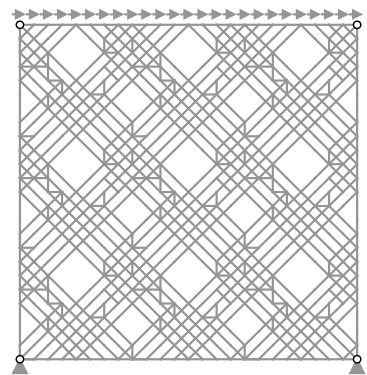

(a) (b)

Figure 11. Near-optimal lattice structure of 8 × 8-grid unit; (a) unit, (b) whole structure

## 6. Conclusion

This paper proposed a methodology to construct an NN model for predicting the responses of periodic lattice structures and apply the pre-trained NN model for small lattice structures to larger structures. We also perform topology optimization of $4 \times 4$-grid units and $8 \times 8$-grid units by the NN-assisted SA algorithm.

The numerical example shows that the filtering process considering the connectivity of adjacent members improves the prediction accuracy of the NN model and that feature selection can compress the NN model while maintaining its performance. In addition, the fact that the NN model has a certain degree of accuracy in applying the pre-trained model to an untrained model indicates that the convolution operation by dividing the data into 4 channels can artificially convert finer lattice patterns into the coarse lattice patterns with close mechanical properties. Furthermore, the trained NN model assists the SA efficiently in searching for solutions with significantly fewer number of structural analyses, compared to SA alone.

## Acknowledgments

This study is partly supported by JSPS KAKENHI No. JP23K04104 and ENEOS scholarship for academic overseas activities.